\documentclass[11pt]{amsart}
\usepackage[dvips]{graphicx} \newcommand{\halmos}{\rule{1ex}{1.4ex}}
\newcommand{\proofbox}{\hspace*{\fill}\mbox{$\halmos$}}
\newdimen\margin   % needed for macros \textdisplay & \ltextdisplay
\def\COMMENT#1{}
\def\TASK#1{}

\def\enddiscard{}
\long\def\discard#1\enddiscard{}

\newcommand{\eps}{\varepsilon}

\newcommand{\prob}{\mathbb{P}}
\newcommand{\ex}{\mathbb{E}}

\newcommand{\G}{G_{n,p}}

\newcommand{\E}{\mathbb{E}}

\newcommand{\U}{{\mathcal U}}
\newcommand{\B}{{\mathcal B}}

\newcommand{\W}{{\mathcal W}}
\newcommand{\C}{{\mathcal C}}
\newcommand{\cP}{{\mathcal P}}

\newtheorem{firsttheorem}{Proposition}

\newtheorem{thm}[firsttheorem]{Theorem}

\newtheorem{corollary}[firsttheorem]{Corollary}

\begin{document}
\title{The order of the largest complete minor in a random graph}
\author{Nikolaos Fountoulakis, Daniela K\"uhn and Deryk
Osthus}
\thanks {N.~Fountoulakis and D.~K\"uhn were supported by the EPSRC, grant no.~EP/D50564X/1} 
\maketitle\vspace{-.8cm}
\begin{abstract}
Let~ccl($G$) denote the order of the largest complete minor in a graph~$G$ (also called 
the contraction clique number)
and let~$G_{n,p}$ denote a random graph on~$n$ vertices with edge probability~$p$.
Bollob\'as, Catlin and Erd\H{o}s~\cite{BCE} asymptotically
determined~ccl($G_{n,p}$) when~$p$ is a constant.  
{\L}uczak, Pittel and Wierman~\cite{LPW} gave bounds on~ccl($G_{n,p}$)
when~$p$ is very close to~$1/n$, i.e.~inside the phase transition. 
We show that for every $\eps>0$ there exists a constant~$C$ such that whenever 
$C/n < p <1-\eps$ then asymptotically almost surely
ccl($G_{n,p}$)$=(1\pm \eps)n /\sqrt{\log_b (np)}$, where
$b:=1/(1-p)$. If $p=C/n$ for a constant $C>1$, then asymptotically almost
surely ccl($G_{n,p}$)$=\Theta(\sqrt{n})$. 
This extends the results in~\cite{BCE} and answers a question
of Krivelevich and Sudakov~\cite{KS}.
\end{abstract}

\section{Introduction}
\subsection{Main results}

A graph~$H$ is a \emph{minor} of~$G$ if for every vertex $h\in H$ there
is a connected subset $B_h\subseteq V(G)$ such that all the~$B_h$ are disjoint
and~$G$ contains an edge between~$B_h$ and~$B_{h'}$ whenever~$hh'$ is an edge of~$H$.
The~$B_h$'s are called the \emph{branch sets}. 
We denote by~ccl$(G)$ the order of the largest complete minor in~$G$.
The study of the largest complete minor contained in 
a given graph has its origins in Hadwiger's conjecture which
states that if the chromatic number of a graph~$G$ is at least~$k$, 
then~$G$ contains $K_k$~minor. It has been proved for $k\le 6$
(see for example~\cite[Chapter~7]{Diest} for a discussion). 

Bollob\'as, Catlin and Erd\H{o}s~\cite{BCE} proved that Hadwiger's conjecture is
true for almost all graphs. For this, they
estimated the typical order of the largest complete minor in 
a graph on~$n$ vertices and compared it with the typical chromatic number of such a graph. 
In particular, they proved that for constant~$p$ and $\eps >0$ 
asymptotically almost surely ccl$(G_{n,p})=(1 \pm \eps)n/\sqrt{\log_b n}$,
where $b:=1/(1-p)$. Here~$G_{n,p}$ is a random graph on~$n$ vertices
where the edges are present independently and with probability~$p$.
We say that an event occurs \emph{asymptotically almost surely} (a.a.s.)
if it occurs with probability tending to~$1$ as~$n$ tends to infinity.

Krivelevich and Sudakov~\cite{KS} considered the order of the largest 
complete minor in a sparser random graph (and more generally in arbitrary
pseudo-random and expanding graphs). 
%They observed that the proof in~\cite{BCE} can be extended to the case
%$p \to 0$ as long as $p$ is not too small, but that it breaks down eventually.
They determined the order of magnitude of~ccl$(G_{n,p})$
as long as~$p\ge n^{\eps-1}$. Our first result determines~ccl$(G_{n,p})$ asymptotically
as long as $p\ge C/n$ and $p=o(1)$.

\begin{thm}\label{thmdense}
For every $\eps>0$ there exists a constant $C=C(\eps)$ such that if $pn\ge C$
and $p=o(1)$, then a.a.s.
$${\rm ccl}(G_{n,p})=(1\pm \eps) \sqrt{\frac{n^2p}{\ln (np)}}. $$ 
\end{thm}

One can combine Theorem~\ref{thmdense} with~\cite{BCE} to obtain
a single formula which allows for constant~$p$ as well. Indeed, let $b:=1/(1-p)$.
If $p=o(1)$ a series expansion gives $\ln b=-\ln (1-p)=p+O(p^2)$. 
Thus 
$$
\sqrt{\frac{n^2p}{\ln (np)}} = \sqrt{\frac{n^2 p}{\ln b\log_b (np)}}
=(1+o(1)){n \over \sqrt{\log_b (np)}}.
$$
Also if~$p$ is constant, then $\log_b n=(1+o(1))\log_b(np)$.   
So altogether we obtain the following.
\begin{corollary}\label{thmdense1}
For every  $\eps>0$ there exists a constant $C=C(\eps)$ such that if $C/n \le p \le 1-\eps$,
then a.a.s. 
$${\rm ccl}(G_{n,p})=(1\pm \eps) \frac{n}{\sqrt{\log_{b} (np)}}. $$ 
\end{corollary}

In the last section of the paper, we estimate~ccl$(G_{n,c/n})$ where $c>1$ is fixed.
Krivelevich and Sudakov~\cite{KS} observed that there are constants~$c_1$ and~$c_2$
such that $c_1\sqrt{n /\log n}\le {\rm ccl}(G_{n,c/n})\le c_2\sqrt{n}$
and asked what the correct order of magnitude is. 

\begin{thm}\label{thm2}
For every $c>1$ there exists a constant $\delta=\delta(c)$ such that a.a.s.
$\delta\sqrt{n}\le {\rm ccl}(G_{n,c/n})\le 2\sqrt{cn}$. 
\end{thm}
Note that the upper bound in Theorem~\ref{thm2} is immediate, since for any graph~$G$, 
the number of edges in any minor of~$G$ is at most~$e(G)$ (see the beginning of Section~\ref{sec:sparse}). 
The same argument shows that the condition that $p\ge c/n$ for some constant $c>1$
is necessary to ensure a complete minor or order~$\Theta(\sqrt{n})$ in~$G_{n,p}$.
This follows from the fact that if $pn\to 1$ the number of edges in any component
is sublinear in~$n$ (see e.g.~\cite[Thm.~6.16]{Bol}). 

%Since the order of the largest complete minor in a graph does not decrease when 
%adding edges, Proposition 1.13 in~\cite{Rg} or Theorem 2.2 in~\cite{Bol}
%imply that the above results carry over to the~$G_{n,M}$ model of random graphs, 
%with $M=p{n \choose 2}$, which is the uniform distribution on all graphs
%on~$n$ vertices and~$M$ edges.%

\subsection{Related results and open questions}

While the influence of the chromatic number on the existence of complete minors is far from clear, 
the corresponding extremal problem for the average degree has been settled for large complete minors:
Thomason~\cite{Thom1} asymptotically determined the smallest average degree~$d(k)$ which guarantees 
the existence of a $K_k$~minor in any graph of average degree at least~$d(k)$.
(The order of magnitude $k \sqrt{\log k}$ of~$d(k)$ was determined earlier in~\cite{Kost,Thom}.)
The extremal graphs are (disjoint copies of) dense random graphs.
Recall that Theorem~\ref{thm2} shows that the behaviour of sparse random graphs
is quite different: in that case~ccl$(G_{n,p})$ has the same order of magnitude as $\sqrt{e(G_{n,p})}$,
a trivial upper bound which holds for any graph. 

There are several results on large complete minors in pseudo-random graphs and expanding graphs.
Thomason~\cite{ATjumbled} introduced the following notion of pseudo-randomness:
a graph~$G$ is \emph{$(p,\beta)$-jumbled} if every induced subgraph~$H$ of~$G$ satisfies
$$\left|e(H)-p\binom{|H|}{2}\right|\le \beta |H|.
$$
He observed that if both~$pn$ and~$(1-p)n$ tend to infinity,
then~$G_{n,p}$ is a.a.s.~$(p,2\sqrt{pn})$-jumbled. 
He also showed%
     \COMMENT{He doesn't actually show this. He just gives a hint how this
can be proved.}
that ccl$(G)\ge (1+o(1))n/\sqrt{\log_b n}$ for every $(p,\beta)$-jumbled graph~$G$
with constant~$p$ and $\beta=O(n^{1-\eps})$. 
A result of Krivelevich and Sudakov~\cite{KS} implies that
ccl$(G)\ge \Omega(\sqrt{n^2p/\log (n\sqrt{p})})$ for $(p,\beta)$-jumbled graphs~$G$
with $\beta=o(np)$. Thus for~$G_{n,p}$ their results only imply
the lower bound in Theorem~\ref{thmdense} up to a multiplicative 
constant if $p \ge n^{\eps-1}$. It would be interesting to know whether
their bound is best possible
or whether (up to a multiplicative constant) the bound in 
Theorem~\ref{thmdense} can be extended to jumbled graphs with appropriate parameters.
Krivelevich and Sudakov~\cite{KS} also considered minors in expanding graphs.
Again, their results only imply the lower bound in Theorem~\ref{thmdense} up to a multiplicative 
constant if $p \ge n^{\eps-1}$. A property closely connected to expansion is that
of having no small separator. Minors in such graphs were considered e.g.~by
Plotkin, Rao and Smith~\cite{PRS}.

%It would be interesting to know whether one can also  find a version
%of Theorem~\ref{thm2} for jumbled graphs.

%Another notion closely connected to pseudo-randomness is that of expansion.
%Krivelevich and Sudakov~\cite{KS} considered the following notion:
%a graph~$G$ on~$n$ vertices is $(t,\alpha)$-expanding if the neighbourhood of every set 
%of $s \le \alpha n /t$ vertices has size at least~$t s$.
%They showed that if $t \ge 10$ and~$\alpha$ is constant, 
%then every $(t,\alpha)$-expanding graph $G$ satisfies
%ccl$(G)\ge \Omega( \sqrt{ n t\log t}/\log n)$.
%They derived the result on minors in random graphs
%mentioned before Theorem~\ref{thmdense} from this (as the expansion~$t$ of
%a random graph is close to its average degree in this range).
%They also showed that if $t \le \log n$, one can improve their bound to 
%$\Omega( \sqrt{ n \log t/\log n})$.
%Similarly to the above question on jumbled graphs, it would be interesting to know whether
%their bound is best possible or whether the bound in Theorem~\ref{thmdense} and~\ref{thm2} can be extended
%to expanding graphs. Note however, that for bounded 
%and slowly growing average degree, the expansion of a random graph is
%$0$ due to the existence of isolated vertices, so one would need a more general notion of expansion
%to get a formulation of our results in this case.

A question that is left open by our results is how the order of the largest complete
minor behaves when~$p$ approaches the
critical point~$1/n$ of the appearance of the giant component. Results in
this direction where proved by \L uczak, Pittel and Wierman~\cite{LPW},
who determined the limiting probability~$g(\lambda)$ that~$G_{n,p}$ is planar when
$np=(1+\lambda n^{-1/3})$. Their results imply that if~$\lambda$ is bounded,
then~$g(\lambda)$ is bounded away from~$0$ and~$1$. (The likely obstacle to planarity
in this case is the existence of a~$K_{3,3}$ minor.)
If $\lambda \to - \infty$, then $g(\lambda) \to 1$.
If $\lambda  \to \infty$, then $g(\lambda) \to 0$.
In fact they proved the stronger result that a.a.s.~ccl$(G_{n,p})$ is unbounded in this case.
%However, it would be very interesting to get a clearer picture of how the largest complete minor grows
%so rapidly as~$np$ increases to~$1+\eps$ for some constant~$\eps.$

Another related problem is that of estimating the order~tcl$(G)$ of the largest
topological clique in a graph~$G$. (Recall that a topological clique 
of order~$k$ is obtained from a~$K_k$ by replacing its edges by internally disjoint paths.)
The typical order of~tcl$(G)$ in a graph~$G$ was first estimated by
Erd\H{o}s and Fajtlowicz~\cite{ErdosF}, who used their bounds to prove that
the conjecture of Haj\'os is false for almost all graphs.
Their estimates were improved by 
Bollob\'as and Catlin~\cite{BC}, who showed that for almost all graphs~tcl$(G)$ 
it is close to~$2\sqrt{n}$. %(see also Exercise 11.9 in~\cite{Bol}). 
Later Ajtai, Koml\'os and Szemer\'edi~\cite{AKS1} proved that as long as the expected 
degree~$(n-1)p$ is at least~$1+\eps$ and is~$o(\sqrt{n})$, 
then a.a.s.~tcl$(G_{n,p})$ is almost as large as the 
maximum degree. Again the likely order of~tcl$(G_{n,p})$ 
when~$p$ is close to~$1/n$ is not known.

Finally, note that our results do not cover the case where $p \to 1$.
Usually, the investigation of~$G_{n,p}$ for such~$p$ is not particularly interesting.
However, any counterexamples to Hadwiger's conjecture are probably rather dense, so in this case it might
be worthwile to investigate the values of the chromatic number and the order of the largest complete minor 
of such a random graph (though it seems rather unlikely that this approach will yield any counterexamples). 

\subsection{Strategy of the proofs}
As in~\cite{BCE}, the upper bound in Theorem~\ref{thmdense} is proved by a first moment argument.
The main difference between the arguments is that in our case, we need to make use of the fact
that the branch sets of a minor have to be connected, whereas this was not necessary
in~\cite{BCE}.

For the lower bound, 
let $k:=n/\sqrt{ \log_b (np)}$ be the function appearing in Corollary~\ref{thmdense1}. 
The proof in~\cite{BCE} for the case when~$p$ is a constant proceeds as follows.
One first shows that a.a.s.~there are~$k$ large pairwise disjoint connected sets~$B_i$
in~$G_{n,p}$. These are used as candidates for the branch sets.
The number~$U_0$ of pairs of~$B_i$ which are not connected by an edge is
then shown to be~$o(k)$. So by discarding a comparatively small number of candidate
branch sets, one can obtain the desired minor. For small~$p$, the main problem is
that~$U_0$ will be much larger than~$k$. However, we can show that~$U_0$ is at most
a small fraction of~$n$. We make use of this as follows. We first find a path~$P$
whose length satisfies $U_0 \ll |P| \ll n$ and which is disjoint from the~$B_i$.
We will divide this path into disjoint subpaths. Our aim is to join most of
those pairs of~$B_i$ which are not yet joined by an edge via one of these subpaths. 
More precisely, we are looking for a matching of size~$U_0-\eps k$
in the auxiliary bipartite random graph~$G^*$
whose vertex classes consist of the unjoined pairs of candidate branch sets and
of the subpaths and where a subpath is adjacent to such an unjoined pair
if it sends an edge to both of the candidate branch sets in this unjoined pair.
There are two difficulties to overcome in order to find such a matching. 
Firstly, some of the~$B_i$ are involved in several unjoined pairs,
so the edges~$G^*$ are not independent. Secondly, if we make the subpaths too short,
then the density of~$G^*$ is not large enough to guarantee
a sufficiently large matching, while if we make the subpaths too long, then
there will not be enough of them. We overcome this by using paths of
very different lengths together with a greedy matching algorithm which
starts off by using short paths to try and join the unjoined pairs. Then in the later stages
the algorithm uses successively longer paths to try and join those pairs
which were not joined in the previous stages until $U_0-\eps k$ of the pairs have been joined.
To ensure that the dependencies between the existence of edges in~$G^*$ are not too large, we 
also remove some of the unjoined pairs from future consideration after each stage
(namely those containing a candidate branch set that is involved in
comparatively many pairs which are still unjoined).

\section{The upper bound of Theorem~\ref{thmdense}}\label{sec:upperthm1}
To prove the upper bound, let 
$$ k := (1+\eps) \sqrt{\frac{n^2p}{\ln (np)}}.$$ 
We use a first moment argument to show that a.a.s. there are no~$k$ 
subsets of~$V$ that can serve as branch sets of a $K_{k}$ minor, where~$V$
denotes the vertex set of~$\G$.%
      \COMMENT{We do need the $+\eps$ term towards the end of the calculation}
Let $B_1,\ldots , B_{k}$ denote any collection of disjoint subsets of 
$V$ and let $s_1, \ldots , s_{k}$ be their sizes. We are not assuming
that the~$B_i$'s cover~$V$ and hence we let~$s$ be the number of
vertices that do not belong to any~$B_i$. So altogether the~$B_i$'s
contain~$n-s$ vertices. We will estimate the probability that 
these sets form the branch sets of a $K_{k}$ minor.

The probability that for all $i <j$, $B_i$ is joined to~$B_j$ by at least one edge is 
$$
\prod_{i<j} \left(1-(1-p)^{s_is_j}\right) \leq \exp \left(-\sum_{i<j}(1-p)^{s_is_j} \right).
$$
An easy argument%
     \COMMENT{For simplicity, let us set $a=(1-p)^{(n-s)^2}$ and let $x_i=s_i/(n-s)$. 
Thus $\sum x_i=1$. Applying the AM-GM inequality to the numbers $a^{x_ix_j}$ 
we obtain $\sum_{i<j} a^{x_ix_j}\geq {k \choose 2}\left(\prod_{i<j}  a^{x_ix_j}\right)^{1/\binom{k}{2}}
={k \choose 2} a^{{k \choose 2}^{-1} \sum_{i<j} x_i x_j}$. 
Now we show that  $\sum_{i<j} x_i x_j\leq {k \choose 2} {1\over k^2}$, applying an induction 
argument on $k$. 
If $k=2$, then $x_2=1-x_1$. But $x_1(1-x_1)$ is maximised for $x_1=1/2$. Assume now that this is 
true when we have $k-1$ variables. We write 
$$ \sum_{i<j} x_i x_j = x_1 (1-x_1) + \sum_{1<i<j}x_i x_j.$$ 
By the induction hypothesis (formally, we apply induction to $x_2/(1-x_1),\dots,x_k/(1-x_1)$
as their sum is 1), the second summand is at most ${k-1 \choose 2} \left(1-x_1 \over k-1 \right)^2$. 
Thus 
$$\sum_{i<j} x_i x_j \leq x_1 (1-x_1) + {k-1 \choose 2} \left(1-x_1 \over k-1 \right)^2.$$
The right-hand side is $x_1 (1-x_1) + {k-2 \over 2(k-1) } (1-x_1)^2 =
{1 \over 2(k-1)}(2(k-1)(x_1-x_1^2)+(k-2)(1-2x_1+x^2_1))=
{1 \over 2(k-1)} (x_1^2((k-2)-2(k-1))+x_1(-2(k-2)+2(k-1))+k-2)=
{1 \over 2(k-1)} (-kx_1^2 +2x_1 + k-2)$. The 
expression in parenthesis is maximized when $x_1=1/k$, yielding 
$$ x_1 (1-x_1) + {k-1 \choose 2} \left(1-x_1 \over k-1 \right)^2 \leq {1\over 2(k-1)} \left( {1\over k} + k-2\right) = 
{1\over 2(k-1)}~{k(k-2)+1\over k} = {k-1 \over 2k} ={k\choose 2}~{1\over k^2}.$$   
}
shows that the sum in the exponent is minimized when the sizes
of the branch sets are as equal as possible. In other words, 
$$ \sum_{i<j} (1-p)^{s_is_j} \geq {k \choose 2}(1-p)^{\left(\frac{n-s}{k}\right)^2}.$$
Therefore,
\begin{equation} \label{ProbJoin} 
\prob (\forall i<j, \; \mbox{$B_i$ is joined to $B_j$}) \leq \exp
\left(-  {k \choose 2}(1-p)^{\left(\frac{n}{k}\right)^2}\right).
\end{equation}
For all $i$, the probability that $B_i$ is connected is at most the expected number of 
spanning trees induced on $B_i$, which is $s_i^{s_i-2}p^{s_i-1}$.
Since all these events are independent for different~$B_i$,
we have
\begin{equation}\label{ProbConn}
\prob \left(\mbox{all $B_i$'s are connected}\right) 
\leq \prod_{i=1}^{k} s_i^{s_i-2} p^{s_i-1} 
\le p^{n-s-k} \prod_{i=1}^{k} s_i^{s_i}.
\end{equation}
Note that for fixed $s_1,\dots,s_k$, there are $\frac{n!}{s!s_1!\cdots s_k!}$ ways of
choosing the candidate branch sets $B_1,\dots,B_k$. 
Since $s_i!\ge (s_i/e)^{s_i}$ by Stirling's formula, we have that
$$
\frac{n!}{s!s_1!\cdots s_k!}=\binom{n}{s}\frac{(n-s)!}{\prod_{i=1}^k s_i!}
\le 2^n \frac{n^{n-s}}{\prod_{i=1}^k \left(\frac{s_i}{e}\right)^{s_i}}
\le \frac{e^{2n} n^{n-s}}{\prod_{i=1}^k s_i^{s_i}}.
$$
Together with~(\ref{ProbJoin}) and~(\ref{ProbConn}) this shows that
for fixed $s_1,\dots,s_k$ the expected number of families of~$k$ disjoint subsets of~$V$
forming the branch sets of a $K_{k}$ minor is at most 
\begin{equation}\label{countBi}
e^{2n} n^{n-s} p^{n-s-k} \exp\left(-  {k \choose 2}(1-p)^{\left(\frac{n}{k}\right)^2}\right)
= e^{2n} (np)^{n-s} \exp\left( -k\ln p- {k \choose 2}(1-p)^{\left(\frac{n}{k}\right)^2}\right).
\end{equation}
But using that $p=o(1)$, we have $1-p \ge e^{-p-p^2} \ge e^{-(1+\eps)p}$
(for the first inequality, see~\cite[p.~5]{Bol}). Together with
our assumption that $pn\ge C=C(\eps)$ this implies that
\begin{align*}
{k \choose 2}(1-p)^{\left(\frac{n}{k}\right)^2} 
\ge \frac{n^2p}{2\ln(np)}e^{-(1+\eps)\ln (np)/(1+\eps)^2} 
= \frac{n}{2\ln(np)} (np)^{1-\frac{1}{1+\eps}}
\ge 3n \ln (np).
\end{align*}
Also, using $p \ge C/n\ge 1/n$, we have 
\begin{align*}
 -k \ln p & \leq 2\ln n \sqrt{\frac{n^2p}{\ln (np)}}=
2\sqrt{n}\ln n \ln (np)\sqrt{\frac{np}{(\ln (np))^3}}\le 
2\sqrt{n}\ln n \ln (np)\sqrt{\frac{n}{(\ln n)^3}}\\
& \le n \ln (np),
\end{align*}
provided that~$C$ is large enough. 
Hence, by~(\ref{countBi}) for fixed $s_1,\dots,s_k$ the expected number of families of~$k$
disjoint subsets of~$V$ forming the branch sets of a $K_{k}$ minor is at most
$$
e^{2n} (np)^{n-s} \exp(n \ln (np)-3n\ln (np)) \le e^{2n} (np)^{-n}.
$$
There are~$\binom{n}{k}$ ways of choosing%
     \COMMENT{Indeed, there are ~$\binom{n}{k}$ ways of choosing a sequence
of~$k$ different numbers $t_1<t_2<\dots<t_k$ from~$[n]$. But setting
$s_1:=t_1$, $s_2:=t_2-t_1,\dots, s_k:=t_k-t_{k-1}$ gives a bijection
between all such sequences and all possible~$s_i$.}
nonnegative integers $s_1,\dots,s_k$ such that $\sum_{i=1}^k s_i\le n$.
Since $e^{2n} (np)^{-n}\binom{n}{k}=o(1)$
this shows that a.a.s.~$G_{n,p}$ does not contain a $K_k$ minor. 

%%%%%%%%%%%%%%%%%%%%%%%%%%%%%%%%%%%%%%%%%%%%%%%%%%%%%%%%%%%%%%%555
%%%%%%%%%%%%%%%%%%%%%%%%%%%%%%%%%%%%%%%%%%%%%%%%%%%%%5

\section{The lower bound of Theorem~\ref{thmdense}}\label{sec:lowerdense2}
We set 
$$ k :=  (1-\eps) \sqrt{\frac{n^2p}{\ln (np)}}, $$
and we will show that there exists a constant $C=C(\eps)$ such that if $np > C$, then
a.a.s.~$G_{n,p}$ contains a $K_{k}$ minor. 

As before, we let~$V$ denote the vertex set of~$G_{n,p}$. 
Let $V'$ be any subset of~$V$ containing $n':=(1-\eps/4)n$ vertices and put
$V'':= V \setminus V'$. First, we construct a family $B_1,\ldots, B_{k'}$ of
$k':=(1+\eps/4)k$ mutually disjoint and connected candidate branch sets in~$G_{n,p}[V']$, 
each having size $t:=(1-\eps/4)n'/k'$. It will turn out that a.a.s.~there are~$k$
of these which can be extended into the branch sets of a~$K_k$ minor.

 \subsection{Constructing $k'$ candidate branch sets}\label{sec:constr}
To construct the $k'$ candidate branch sets $B_1,\dots,B_{k'}$ we will choose
a long path in $G_{n,p}[V']$ and divide it into consecutive subpaths,
each consisting of~$t$ vertices. We will obtain this path via a
two-round exposure of $G_{n,p}[V']$. More precisely, let $p',p''$ be such
that $p=p' + p'' - p'p''$. In our two-round exposure of 
$G_{n,p}[V']$ we first generate a random graph $G_{n,p'}[V']$ and
then independently $G_{n,p''}[V']$. We then take their union, ignoring multiple edges
(i.e. regarding them as a single edge). 

In our case, we will take $p'$ to be a small proportion of~$p$, whereas~$p''$
will be almost equal to~$p$. However, we have to make sure that~$p'$ is not too small
either, so that $G_{n,p'}[V']$ contains a long path.
Ajtai, Koml\'os and Szemer\'edi~\cite{AKS} and de la Vega~\cite{Vega}
independently proved that there exists a constant $C'=C'(\eps)$ such that if $p'|V'|=C'$
then $G_{n,p'}[V']$ contains a path~$P'$ on $(1-\eps/4)|V'|$ vertices
(see also Theorem 8.1 in~\cite{Bol}).  
It is easily seen that if we choose~$C$ sufficiently large compared to~$C'$ and~$1/\eps$
then
\begin{equation}\label{eq:p''}
p''=\frac{p-p'}{1-p'}\ge p-p'\ge (1-\eps^2)p.
\end{equation}  
We divide the path~$P'$ into consecutive disjoint subpaths $B_1,\ldots, B_{k'}$, each containing~$t$
vertices. We shall use these~$B_i$'s as candidate branch sets. In what follows, we assume
that~$C$ was chosen to be sufficiently large compared to~$1/\eps$ and~$n$
was chosen to be sufficiently large compared to~$C$ for our estimates to hold.  

\subsection{Estimating the number of unjoined pairs of candidate branch sets}\label{sec:est}
Our aim now is to show that~$k$ of the candidate branch sets $B_1,\dots,B_{k'}$
are joined pairwise to each other (by edges or paths through~$V''$) and thus form a $K_{k}$~minor. 
First, we will estimate the number of the pairs $B_i,B_j$ that are not joined
by an edge in the second round of the two-round exposure of~$G_{n,p}[V']$. So let~$\U_0$ denote the set
of all those pairs~$B_i,B_j$ which are not joined by an edge in~$G_{n,p''}[V']$
and let $U_0:=|\U_0|$.
Consider the auxiliary graph~$G'_0$ whose vertices are the~$B_i$'s and whose edges
correspond to the pairs in~$\U_0$.
So~$G'_0$ is a binomial random graph on~$k'$ vertices whose edges occur independently
with probability $q:=(1-p'')^{t^2}$. Note that if~$\eps$ is small enough then%
     \COMMENT{Indeed, to see the 2nd inequality note that
$(1+\eps/4)^2(1-\eps)=(1+\eps/2+\eps^2/16)(1-\eps)=1-\eps/2-\eps^2/2+\eps^2/16-\eps^3/16
\le 1-\eps/2$.}
\begin{align}\label{lowerboundt}
t^2 & = {(1-\eps/4)^4 \over (1+\eps/4)^2 (1-\eps)^2} {\ln (np) \over p}
\ge \left( \frac{1-\eps/2}{(1+\eps/4)(1-\eps)}\right)^2 {\ln (np) \over p}
\ge (1+\eps/4)^2 \frac{\ln (np)}{p} \nonumber\\
& \ge (1+\eps/2)\frac{\ln (np)}{p}
\end{align}
and so
\begin{equation*}
q=(1-p'')^{t^2}\le  e^{-p'' t^2} \stackrel{(\ref{eq:p''})}{\le} e^{-(1-\eps^2)pt^2} \le 
\left( \frac{1}{np}\right)^{1+\eps/3}.
\end{equation*}
Thus
\begin{eqnarray}\label{eq:exU0}
\E (U_0)={k' \choose 2} q  \le { n^2 p \over \ln (np) } \frac{1}{(np)^{1+\eps/3}}
\le {n \over (np)^{\eps/3}}.
\end{eqnarray}
Suppose first that $\E(U_0)\le n^{2/5}$. Then Chebyshev's inequality implies%
    \COMMENT{Chebyshev's inequality implies that a $Bin(n,p)$ random variable $X$
satisfies that $\prob(|X-\E(X)|\ge \eps \E(X))
\le \frac{Var(X)}{\eps^2(\E(X))^2}=\frac{npq}{\eps^2 n^2p^2}\le \frac{1}{\eps^2\E(X)}$.
So we can apply this with $\eps:=n^{2/5}/\E(U_0)$.} 
that a.a.s.~$U_0\le 2n^{2/5}\le \eps^2 k'$. Thus by deleting one candidate
branch set from every pair in~$\U_0$ we obtain a $K_k$ minor.
So we may assume that~$\E(U_0)\ge n^{2/5}$. (In fact, a straightforward but tedious
calculation shows that this is always the case.)
In particular, $\E (U_0) \rightarrow \infty$ as $n \rightarrow \infty$ and
so a.a.s.
\begin{equation}\label{eq:U0}
U_0 =(1\pm \eps) \E (U_0)
\end{equation}
by Chebyshev's inequality. 
Note that in the case when $p\ge 1/(\ln n)^2$ inequalities~(\ref{eq:exU0})
and~(\ref{eq:U0}) together imply that a.a.s
$$ \eps^2 k'\ge \frac{\eps^2 n}{2}\sqrt{\frac{p}{\ln n}}\ge
\frac{\eps^2 n}{(\ln n)^2}\ge 2n^{1-\eps/3}(\ln n)^{2\eps/3}
\ge \frac{2n^{1-\eps/3}}{p^{\eps/3}}\ge U_0.
$$
Thus by deleting one candidate branch set from every pair in~$\U_0$
we obtain a $K_k$ minor.  

So in what follows, we may assume that $p\le 1/(\ln n)^2$.
We will show that most of the unjoined pairs of candidate branch sets (i.e.~most of
the pairs in~$\U_0$) can be joined using~$V''$. More precisely, 
for every pair $u\in \U_0$ we will try to find a path in~$G_{n,p}[V'']$
which sends an edge to both candidate branch sets belonging to~$u$. 
Of course these paths have to be disjoint for different pairs. So they will be
chosen as follows: we consider a longest path in $G_{n,p}[V'']$,
divide it into consecutive disjoint subpaths of various lengths and use them to successively 
join the pairs of candidate branch sets in~$\U_0$.
Here the fact that by~(\ref{eq:exU0}) and~(\ref{eq:U0}) a.a.s.~the number
$U_0=|\U_0|=e(G'_0)$ of unjoined pairs is much smaller than~$|V''|$ is crucial.

More formally, we proceed as follows. Let%
    \COMMENT{(\ref{eq:U0}) together with our assumption that
$\E(U_0)\ge n^{2/5}$ implies that $i^*\to \infty$.}
$$i^*:=\frac{\ln U_0-\ln(n^{1/3})}{\ln 8}\stackrel{(\ref{eq:U0})}{\le} \ln n.
$$
For each $i=1,\dots,i^*$ let
\begin{equation}\label{eq:Ui}
U_i:=U_0/8^i.
\end{equation}
Thus $U_{i^*}=n^{1/3}$. The results in~\cite{AKS, Vega} imply that~$G_{n,p}[V'']$
contains a path on $|V''|/2\ge \eps^2 n$ vertices.
We divide this path into disjoint consecutive subpaths $Q_1, Q_2, \ldots,Q_{i^*}$ such that%
     \COMMENT{$|Q_i|\to \infty$ as $|Q_i|\ge |Q_{i^*}|\ge \frac{\eps^2 n}{2^{\ln U_0/\ln 8}\ln(np)}
= \frac{\eps^2 n}{2^{\log U_0/\log 8}\ln(np)}=\frac{\eps^2 n}{U_0^{1/3}\ln(np)}
\ge \frac{\eps^2 n^{2/3}}{\ln(np)}$ by~(\ref{eq:exU0}).}
\begin{equation}\label{eq:Qi}
|Q_i|=\frac{\eps^2 n}{2^{i}\ln(np)} \ \ \ \ \ \text{and so}\ \ \ \ \
|Q_i|=\frac{|Q_1|}{2^{i-1}}.
\end{equation}
For each~$i=1,\dots,i^*$ we divide~$Q_i$ further into a set~$\cP_i$ of~$U_{i-1}$
disjoint smaller consecutive paths each containing $\ell_i:=|Q_i|/U_{i-1}$
vertices.%
     \COMMENT{Note that $\ell_i\ge
\frac{\eps^2n 4^{i-1}}{2\ln(np)U_0}\ge
\frac{\eps^2n 4^{i-1}}{2\ln(np)2\ex(U_0)}\ge \frac{\eps^2 n 4^{i-1}(np)^{\eps/3}}{4\ln(np)n}$,
which is at least a large constant. So the def of $\ell_i$ makes sense.}
Thus
\begin{equation}\label{eq:li}
\ell_i=\frac{|Q_i|}{U_{i-1}}=\frac{8^{i-1}|Q_1|}{2^{i-1}U_0}=\frac{4^{i-1}|Q_1|}{U_0}
\ \ \ \ \ \text{and so}\ \ \ \ \ \ell_i=4^{i-1}\ell_1.
\end{equation}
Note that~$\ell_i$ and~$|Q_i|$ are large. So in particular viewing them as
integers does not affect the calculations.
Given a path~$P$ in~$G_{n,p}[V'']$ and a candidate branch set~$B$ we write
$P\sim B$ if some vertex in~$B$ sends an edge to some vertex on~$P$.
Given a pair $u\in\U_0$ consisting of candidate branch sets~$A$ and~$B$,
we write $P\sim u$ if $P\sim A$ and $P\sim B$.

Roughly speaking, for each $i=1,2,\ldots,i^*$ in turn, we will use a greedy algorithm to join all
but precisely~$U_i$ of all those pairs in~$\U_0$ that are still unjoined, using
one path in~$\cP_i$ for each such pair $u\in\U_0$, i.e.~we want to find a path $P\in\cP_i$
such that $P\sim u$ (see Sections~\ref{sec:greedy1} and~\ref{sec:ana}
for details). We call each iteration of our algorithm a \emph{stage}. Thus after the $i$th
stage a.a.s.~precisely~$U_i$ pairs in~$\U_0$ will still be unjoined and we aim to join
precisely $7U_i/8=U_i-U_{i+1}$ of these during the next stage. Since~$U_{i^*}=n^{1/3}$
this means that after~$i^*$ stages we obtain a complete minor of order
$k'-n^{1/3}\ge k$ by removing one candidate branch set
for every of the~$U_{i^*}$ pairs in~$\U_0$ which are still unjoined.

\subsection{A greedy algorithm for joining the candidate branch sets}\label{sec:greedy1} 
Let~$\hat{p}$ denote the probability
that a given vertex $v\in V''$ is joined to a given candidate branch set.
Thus $\hat{p}=1-(1-p)^t$. But%
     \COMMENT{The 1st inequality holds by induction since
$(1-p)^{t+1}=(1-p)(1-p)^t\ge (1-p)(1-t p)\ge 1-(t+1)p$.
The 2nd inequality below can by proved by induction:
$(1-p)^{t+1}=(1-p)^t(1-p)\le (1-pt+(pt)^2/2)(1-p)=1-pt+(pt)^2/2-p+p^2t-p^3t^2/2
\le 1-p(t+1)+p^2(t^2+2t+1)/2$.}
\begin{equation}\label{eq:estp}
1-pt\le (1-p)^t\le 1-pt+(pt)^2/2.
\end{equation}
Similarly as in~(\ref{lowerboundt}) one can show that $pt^2=O(\ln (np))$.
So our assumption that $p\le 1/(\ln n)^2$ implies~$pt=o(1)$. Thus 
\begin{equation}\label{eq:boundhatp}
\hat{p} =(1+o(1))pt.
\end{equation}
Suppose that $1\le i< i^*$ and after the $(i-1)$th stage of our algorithm we
have a set $\U_{i-1}\subseteq \U_0$ of unjoined pairs such that $|\U_{i-1}|=U_{i-1}$ and that we
now wish to run the $i$th stage of our algorithm.
(If $i=1$ then $\U_0$ is the set defined at the beginning of Section~\ref{sec:est}.)
So we now wish to use the
paths in~$\cP_{i}$ to join precisely $7U_{i-1}/8$ of the pairs in~$\U_{i-1}$.
(The set of remaining pairs in~$\U_{i-1}$ will then be~$\U_i$.)
However, before running the actual algorithm we will first discard all those
pairs in~$\U_{i-1}$ which contain a candidate branch set lying in too many
pairs from~$\U_{i-1}$. This will reduce the dependencies between the successes of the
individual steps of the algorithm within a stage. More precisely, let~$G'_{i-1}$ denote the
spanning subgraph of~$G'_0$ whose edge set corresponds
to~$\U_{i-1}$. ($G'_0$ is the graph defined in Section~\ref{sec:est}.)
For all $i=1,\dots,i^*$ put
\begin{equation}\label{eq:Deltai}
\Delta_{i-1}:=\frac{8\E (U_0)}{\eps^{3/2} 4^{i-1}k'}.
\end{equation}
Let~$\B_{i-1}$ be the set of all those candidate branch sets whose degree
in~$G'_{i-1}$ is greater than~$\Delta_{i-1}$. In other words, $\B_{i-1}$
consists of all those candidate branch sets which lie in more
than~$\Delta_{i-1}$ pairs from~$\U_{i-1}$.
As $|\B_{i-1}|\Delta_{i-1}\le 2e(G'_{i-1})=2U_{i-1}$, we have that
\begin{equation*}
\frac{|\B_{i-1}|}{k'}\stackrel{(\ref{eq:Deltai})}{\le}
\frac{\eps^{3/2} 4^{i-1}U_{i-1}}{4 \E (U_0)}
\stackrel{(\ref{eq:U0})}{\le} \frac{\eps^{3/2} 4^{i-1}U_{i-1}}{2 U_0}
\stackrel{(\ref{eq:Ui})}{=} \frac{\eps^{3/2}}{ 2^{i}}.
\end{equation*}
Thus
\begin{equation}\label{eq:sumBi}
\sum_{j\ge 0} |\B_{j}|\le \eps^{3/2} k'\sum_{j=0}^\infty 2^{-(j+1)}= \eps^{3/2} k'.
\end{equation}
In the $i$th stage of our algorithm, we discard all those pairs in~$\U_{i-1}$
which contain at least one candidate branch set from~$\B_{i-1}$, i.e.~we
will consider all these pairs as being joined (even if they are not joined).
After running all the~$i^*$ stages of our algorithm we will make up for this by deleting all
the candidate branch sets in $\bigcup_{j\ge 0}\B_j$. By~(\ref{eq:sumBi})
we do not loose too many of the candidate branch sets in this way.
More formally, in the $i$th stage of our algorithm we proceed as follows:
\begin{itemize}
\item[{\rm (i)}] Consider all those pairs in~$\U_{i-1}$ which contain at least one
candidate branch set from~$\B_{i-1}$. If there are at least $7U_{i-1}/8$
such pairs, delete precisely $7U_{i-1}/8$ of them and stop the $i$th
stage of the algorithm.
\item [{\rm (ii)}] If there are less than $7U_{i-1}/8$
such pairs, let $\U^*_{i-1}$ denote the set of all these
pairs, let $\U'_{i-1}:=\U_{i-1}\setminus \U^*_{i-1}$
and $U'_{i-1}:=|\U'_{i-1}|$.
So we are aiming to join precisely~$7U_{i-1}/8-|\U^*_{i-1}|=:U''_{i-1}$
pairs in~$\U'_{i-1}$ by using paths in~$\cP_i$. To do this, we proceed as follows:
\begin{itemize}
\item[{\rm (a)}]  Let $u_1,\ldots, u_ {U'_{i-1}}$ be an ordering of the pairs
in~$\U'_{i-1}$ and let $P_1, \ldots, P_{U_{i-1}}$ be an ordering of the
paths in~$\cP_i$. 
\item[{\rm (b)}] For each~$P_j$ in turn, consider the set of all remaining pairs
in~$\U'_{i-1}$. If there are precisely~$U'_{i-1}-U''_{i-1}$ such pairs, stop
the $i$th stage of the algorithm. (This means that we have successfully joined the desired
number~$U''_{i-1}$ of pairs in~$\U'_{i-1}$.) If there are more than~$U'_{i-1}-U''_{i-1}$ such pairs,
then among all these pairs let~$u_\ell$ be the one with the smallest index such that
$P_j\sim u_\ell$ (if such an~$u_\ell$ exists), remove~$u_\ell$ and continue
with~$P_{j+1}$. If there is no such~$u_\ell$, then continue with~$P_{j+1}$.
In the latter case we say that~$P_j$ \emph{fails}.
\end{itemize}
\end{itemize}
$\U_i$ will be the subset of~$\U_{i-1}$ obtained after the $i$th stage of
our algorithm. So our aim is to show that a.a.s.~$|\U_{i}|=U_i$ for all $i=1,\dots,i^*$.
To do this, it suffices to show that, for all $i=1,\dots,i^*$,
if we are in Part~(ii) of the algorithm, then in
Part~(b) a.a.s.~to~$U''_{i-1}$ paths in~$\cP_i$ we assign pairs in~$\U'_{i-1}$.
As $|\cP_i|-U''_{i-1}=U_{i-1}-U''_{i-1}\ge U_{i-1}/8$ we have that
$$
\prob(\ge |\cP_i|-U''_{i-1} \text{ paths in } \cP_i \text{ failed})\le
\prob(\ge U_{i-1}/8 \text{ paths in } \cP_i \text{ failed})
$$
and so 
it suffices to show that the latter probability is sufficiently small.
Note that this probability is determined by the bipartite subgraph~$G^*$ of~$G_{n,p}$
whose first vertex class is~$V'$ and whose second vertex class
consists of all the vertices on the paths in~$\cP_i$. Moreover, note
that none of the events considered in Sections~\ref{sec:constr} and~\ref{sec:est}
affects the existence of edges in~$G^*$. Furthermore, the success of different
stages of the algorithm depends on edge-disjoint subgraphs of~$G^*$.
So by the principle of deferred decisions, when analysing the success
of the $i$th stage, we may assume that~$\U_{i-1}$ is fixed (and thus also~$\U'_{i-1}$)
and that the bipartite graph~$G^*_i$ between the vertices
in the candidate branch sets belonging to~$\U_{i-1}$ and the vertices in the paths
in~$\cP_i$ is a binomial random graph with edge probability~$p$. Moreover,
when analysing the success of the $i$th stage the events we consider can be
viewed as sets of bipartite graphs having the same vertex classes as~$G^*_i$.

\subsection{Analysis of the greedy algorithm}\label{sec:ana}
Recall that all the paths in~$\cP_i$ consist of~$\ell_i$ vertices. 
Let~$\hat{p}_i$ denote the probability that a given
path~$P\in \cP_i$ is joined to a given candidate branch set~$B$, i.e.~that
$P \sim B$. Thus
\begin{equation}\label{eq:upperphat1}
\hat{p}_i=1-(1-\hat{p})^{\ell_i}\le 1-\left(1-\ell_i\hat{p}\right)=\ell_i\hat{p}.
\end{equation}
Moreover, if $\hat{p}\ell_i\le 1$ then as in~(\ref{eq:estp})
\begin{equation}\label{eq:lowerphat1}
\hat{p}_i=1-(1-\hat{p})^{\ell_i}\ge \ell_i\hat{p}-\frac{(\ell_i\hat{p})^2}{2}
\ge \frac{\ell_i\hat{p}}{2}
\end{equation}
while if $\hat{p}\ell_i\ge 1$ then%
     \COMMENT{works if $\ell_i\ge 2$}
\begin{equation}\label{eq:lowerphat2}
\hat{p}_i=1-(1-\hat{p})^{\ell_i}\ge 1-(1-1/\ell_i)^{\ell_i}\ge 1/2.
\end{equation}
As $|\cP_i|=U_{i-1}$ we have that
\begin{align}\label{eq:failedP1}
\prob(\ge U_{i-1}/8 \text{ paths in } \cP_i \text{ failed}) &=
\sum_{A\subseteq [U_{i-1}],\ |A|\ge U_{i-1}/8} \prob(P_j\text{ fails iff } j\in A)\nonumber\\
& \le 2^{U_{i-1}} \max_{A\subseteq [U_{i-1}],\ |A|\ge U_{i-1}/8} \prob(P_j\text{ fails iff } j\in A).
\end{align}
So let us now consider any $A\subseteq [U_{i-1}]$ with $|A|\ge U_{i-1}/8$.
For all $j\in A$ put $A_j:=A\cap [j-1]$. Then
\begin{align}\label{eq:failedP1A}
\prob(P_j\text{ fails iff } j\in A)\le \prob(P_j\text{ fails } \forall j\in A)=
\prod_{j\in A} \prob(P_j\text{ fails}\mid P_r \text{ fails }\forall r\in A_j).
\end{align}
To estimate the latter probability, we split up the event that~$P_r$ fails for all
$r\in A_j$ further. Let~$\W_j$ be the set of all bipartite graphs whose first vertex
class~$V^*$ consists of all the vertices in the candidate branch sets belonging to~$\U_{i-1}$
and whose second vertex class is~$V(P_1)\cup\dots\cup V(P_{j-1})$.
We view each $W\in \W_j$ as the event that the bipartite subgraph of~$G^*_i$ induced
by~$V^*$ and $V(P_1)\cup\dots\cup V(P_{j-1})$ equals~$W$ ($G^*_i$ was defined at the end of
Section~\ref{sec:greedy1}). Note that
for every $W\in\W_j$ the intersection of the event~$W$ and the event that~$P_r$
fails for all $r\in A_j$ is either~$W$ or empty. Let $\W^*_j$ be the set of all those
$W\in \W_j$ for which the former holds. Then
\begin{align}\label{eq:failedcondPj}
\prob(P_j\text{ fails} \mid P_r & \text{ fails } \forall r\in A_j) =\nonumber \\
& =\sum_{W\in \W^*_j} \prob(P_j\text{ fails}\mid (P_r \text{ fails }\forall r\in A_j)\cap W)
\prob(W\mid P_r \text{ fails }\forall r\in A_j)\nonumber\\
& \le \left(\max_{W\in \W^*_j} \prob(P_j\text{ fails}\mid (P_r \text{ fails }\forall r\in A_j)\cap W)\right)
\sum_{W\in \W^*_j} \prob(W\mid P_r \text{ fails }\forall r\in A_j)\nonumber\\
& = \max_{W\in \W^*_j} \prob(P_j\text{ fails}\mid W).
\end{align}
Now note that for each $W\in \W_j^*$ there is
a set~$\U_W\subseteq \U'_{i-1}$ such that whenever the event~$W$ occurs, the set of remaining
pairs in~$\U'_{i-1}$ after we have considered $P_1,\dots,P_{j-1}$ in the $i$th stage of
our algorithm is precisely~$\U_W$ (i.e.~the set of remaining pairs in~$\U'_{i-1}$
is the same for every bipartite graph belonging to the event~$W$). Thus
\begin{align}\label{eq:failedcondW}
\prob(P_j\text{ fails }\mid W) & =\frac{\prob( P_j\text{ fails } \cap W)}{\prob( W )}
=\frac{\prob( (P_j\not\sim u \;\forall u \in \U_W)\cap W)}{\prob( W )}\nonumber\\
& =\frac{\prob( P_j\not\sim u \;\forall u \in \U_W) \prob(  W)}{\prob( W )}
=\prob( P_j\not\sim u \;\forall u \in \U_W).
\end{align}
To estimate the latter probability, we may assume without loss of generality that~$\U_W$
consists of $u_1,\dots,u_{|\U_W|}$. Note that
\begin{equation}\label{eq:UW}
|\U_W|> U_{i-1}/8
\end{equation}
as otherwise we would have stopped the $i$th stage of our algorithm.
Moreover,
\begin{align}\label{eq:prob}
\prob( P_j\not\sim u \;\forall u \in \U_W)
& =\prod_{r=1}^{|\U_W|} \prob( P_j\not\sim u_r \mid P_j\not\sim u_s \ \forall s< r)\nonumber\\
& = \prod_{r=1}^{|\U_W|} \left(1-\prob( P_j\sim u_r \mid P_j\not\sim u_s \ \forall s< r)\right)\nonumber\\
& = \prod_{r=1}^{|\U_W|} \left(1-\frac{\prob( P_j\sim u_r,\; P_j\not\sim u_s \ \forall s< r)}
{\prob(P_j\not\sim u_s \ \forall s< r)}\right).
\end{align}
Note that the definition of~$\hat{p}_i$ implies that
$$\prob( P_j\sim u_1)=\hat{p_i}^2.
$$
So if all of the above events were independent, then the value of each bracket would be equal
to~$1-\hat{p}_i^2$. Our next aim is to show that this is approximately true
(see~(\ref{eq:prob3}) below). Hence we now consider any $2\le r\le {|\U_W|}$.
Let~$B$ and~$B'$ denote the candidate branch sets
belonging to~$u_r$. Let~$\U_{r,1}$ be the set of all those~$u_s$ with $s<r$ which contain~$B$
or~$B'$. Let $\U_{r,2}$ consist of those~$u_s$ with $s<r$ not in $\U_{r,1}$
for which at least one of their candidate branch sets belongs to a pair in $\U_{r,1}$.
Finally, we let $\U_{r,3}:=\{u_1,\dots,u_{r-1}\}\setminus (\U_{r,1}\cup \U_{r,2})$
be the set of the remaining pairs. Note that
\begin{eqnarray}\label{eq:prob0}
\prob(P_j\sim u_r,\; P_j\not\sim u_s \ \forall s< r)&=& 
\prob(P_j\sim u_r,\; P_j\not\sim u \ \forall u\in \U_{r,1}\cup \U_{r,2}\cup \U_{r,3})\nonumber\\
&=& \prob(P_j\sim u_r,\; P_j\not\sim u \ \forall u\in \U_{r,1}\cup \U_{r,3})\nonumber\\
&=& \prob(P_j\sim u_r,\; P_j\not\sim u \ \forall u\in \U_{r,1})
\prob(P_j\not\sim u \ \forall u\in \U_{r,3}). 
\end{eqnarray}
Indeed, to see the second equality let~$\C_r$ denote the set of all those
candidate branch sets~$C$ with $(B,C) \in\U_{r,1}$ or $(B',C) \in\U_{r,1}$.
Recall that $P_j \sim u_r$ means $P_j\sim B$ and $P_j\sim B'$. So if $P_j \sim u_r$
and $P_j\not\sim u$ for all $u\in \U_{r,1}$ then
$P_j\not\sim C$ for all $C\in\C_r$. But this means that automatically
$P_j\not \sim u$ for any $u \in \U_{r,2}$. The last equality holds since the events
$P_j \sim u_r, \; P_j \not \sim u \; \forall u\in \U_{r,1}$
and $P_j \not \sim u \; \forall u\in \U_{r,3}$ are 
independent as they involve disjoint families of candidate branch sets.
Now
\begin{equation}\label{eq:prob1}
\prob(P_j\sim u_r,\; P_j\not\sim u \ \forall u\in \U_{r,1})=
\prob(P_j\sim B,\; P_j\sim B',\; P_j\not\sim C \ \forall C\in \C_r)=
\hat{p}_i^2(1-\hat{p}_i)^{|\C_r|}.
\end{equation}
We wish to show that the right hand side of~(\ref{eq:prob1}) is at least
$\hat{p}^2_i/2$. So we have to show that $\hat{p}_i|\C_r|$ is small.
To do this, note that $B,B'\notin \B_{i-1}$ since
$u_r=(B,B')\in \U_W\subseteq \U'_{i-1}$. Thus
$|\C_r|\le d_{G'_{i-1}}(B)+d_{G'_{i-1}}(B')\le 2\Delta_{i-1}$ and so
\begin{eqnarray}\label{eq:hatpiCr}
\hat{p}_i|\C_r| &\le &  2\hat{p}_i \Delta_{i-1}
\stackrel{(\ref{eq:upperphat1}), (\ref{eq:Deltai})}{\le}
2\hat{p}\ell_i \frac{8\E (U_0)}{\eps^{3/2} 4^{i-1}k'}
\stackrel{(\ref{eq:boundhatp}), (\ref{eq:li})}{\le}
3pt \frac{|Q_1|}{U_{0}}\cdot \frac{8\E (U_0)}{\eps^{3/2}k}\nonumber\\
& \stackrel{(\ref{eq:Qi}),(\ref{eq:U0})}{\le} &
\frac{3pn}{k}\cdot \frac{\eps^2n}{2\ln(np)}\cdot
\frac{10}{\eps^{3/2}k} \le 20\sqrt{\eps}.
\end{eqnarray}
Together with~(\ref{eq:prob1}) and the fact that
$1-\hat{p}_i\ge e^{-\hat{p}_i-\hat{p}_i^2}\ge e^{-2\hat{p}_i}$
(see e.g.~\cite[p.~5]{Bol}) this implies that
\begin{equation}\label{eq:prob2}
\prob(P_j\sim u_r,\; P_j\not\sim u \ \forall u\in \U_{r,1})\ge
\hat{p}_i^2 e^{-2\hat{p}_i|\C_r|}\ge \hat{p}_i^2/2.
\end{equation}
Thus
\begin{eqnarray}\label{eq:prob3}
\frac{\prob( P_j\sim u_r,\; P_j\not \sim u_s \; \forall s< r)}
{\prob(P_j\not\sim u_s \; \forall s< r)}&= &
\frac{\prob( P_j\sim u_r,\; P_j\not \sim u_s \; \forall s< r)}
{\prob(P_j\not\sim u \ \forall u\in \U_{r,1}\cup \U_{r,2}\cup \U_{r,3})}\nonumber \\
& \ge & \frac{\prob( P_j\sim u_r,\; P_j\not \sim u_s \; \forall s< r)}
{\prob(P_j\not\sim u \ \forall u\in \U_{r,3})}\nonumber \\
&\stackrel{(\ref{eq:prob0})}{=}&
\prob(P_j\sim u_r,\; P_j\not\sim u \ \forall u\in \U_{r,1})
\stackrel{(\ref{eq:prob2})}{\ge} \hat{p}^2_i/2
\end{eqnarray}
and hence
$$\prob( P_j\not\sim u\ \forall u \in \U_W)\stackrel{(\ref{eq:prob}),(\ref{eq:prob3})}{\le}
\left(1-\hat{p}^2_i/2\right)^{|\U_W|}\stackrel{(\ref{eq:UW})}{\le}
\exp\left(-\hat{p}^2_iU_{i-1}/16\right).
$$
Together with~(\ref{eq:failedP1A})--(\ref{eq:failedcondW}) this in turn implies that 
\begin{align*}
\prob(P_j\text{ fails iff } j\in A)\le \exp\left( -\hat{p}^2_iU_{i-1}|A|/16\right)
\le \exp\left(-\hat{p}^2_i (U_{i-1})^2/2^7\right)
\end{align*}
for any subset $A\subseteq [U_{i-1}]$ with $|A|\ge U_{i-1}/8$.
Together with~(\ref{eq:failedP1}) this gives
\begin{align}\label{eq:failedP12}
\prob(\ge U_{i-1}/8 \text{ paths in } \cP_i \text{ failed}) &\le
2^{U_{i-1}} \exp\left(-\hat{p}^2_i(U_{i-1})^2/2^{7}\right)\le
\exp\left(U_{i-1}(1-\hat{p}^2_iU_{i-1}/2^{7})\right).
\end{align}
Together with~(\ref{eq:lowerphat2}) this implies that 
$\prob(\ge U_{i-1}/8 \text{ paths in } \cP_i \text{ failed})\le e^{-U_{i-1}}$
if $\hat{p}\ell_i\ge 1$. So we may assume that $\hat{p}\ell_i\le 1$
and thus $\hat{p}_i\ge \ell_i\hat{p}/2$ by~(\ref{eq:lowerphat1}).
But then
\begin{eqnarray}\label{eq:hatpiUi}
\hat{p}_i^2U_{i-1} &\ge & \left(\frac{\ell_i\hat{p}}{2}\right)^2 U_{i-1}
\stackrel{(\ref{eq:li})}{=}\frac{4^{i-1} \ell_1|Q_i|\hat{p}^2}{4}
\stackrel{(\ref{eq:boundhatp})}{\ge}
\frac{4^{i-1} |Q_i|p^2t^2}{5}
\stackrel{(\ref{eq:Qi})}{\ge}
\frac{\eps^2 n}{2\ln(np)} \cdot\frac{p^2n^2}{6 k^2}\nonumber\\
& \ge & \frac{\eps^2 pn}{12}\ge 2^{8}.
\end{eqnarray}
Together with~(\ref{eq:failedP12}) this implies that in both cases we have 
\begin{equation}\label{eq:failedP13}
\prob(\ge U_{i-1}/8 \text{ paths in } \cP_i \text{ failed})\le e^{-U_{i-1}}=o(1/\ln n).
\end{equation}
The final bound follows from the fact that
$U_{i-1}\ge U_{i^*}=n^{1/3}$. Note that the derivation of~(\ref{eq:hatpiUi}) made
essential use of the fact that the probability~$\hat{p}_i^2$ that a given
path joins a given pair of candidate branch sets grows quadratically in~$\ell_i$
in each stage. This is the reason why it is worthwhile to use fewer but longer
paths in each stage.

As $i^*\le \ln n$, (\ref{eq:failedP13}) implies that
a.a.s.~for all $i=1,\dots,i^*$ our algorithm produces sets~$\U_i$ with $|\U_i|=U_i$.
Recall that each pair of candidate branch sets in~$\U_0$ which has not been joined
by the algorithm either lies in~$\U_{i^*}$ or contains
some candidate branch set belonging to $\bigcup_{j\ge 0}\B_j$.
So by deleting all the candidate branch sets in $\bigcup_{j\ge 0}\B_j$ as well as one
candidate branch set from every pair in~$\U_{i^*}$ we obtain
a complete graph of order $k'-U_{i^*}-|\bigcup_{j\ge 0}\B_j|\ge k$ as minor.
(The last inequality follows from~(\ref{eq:sumBi}) and the fact that $U_{i^*}=n^{1/3}$.) 

%%%%%%%%%%%%%%%%%%%%%%%%%%%%%%%%%%%%%%%%%%%%%%%%%%%%%%%%%%%%%%%%%%

\section{Proof of Theorem~\ref{thm2}} \label{sec:sparse}

To verify the upper bound of Theorem~\ref{thm2} note that 
a.a.s.~$G_{n,c/n}$ has at most~$cn$ edges. But this implies that a.a.s.~the order~$k$
of the largest complete minor in~$G_{n,c/n}$ satisfies $\binom{k}{2}\le cn$,
i.e.~$k\le 2\sqrt{cn}$.

To prove the lower bound of Theorem~\ref{thm2} we are
aiming to show that for every constant $c>1$ there exists a constant $\delta = \delta (c)>0$
such that a.a.s. $G_{n,c/n}$ contains a $K_{k}$ minor where
$$k:=\delta\sqrt{n}.$$
(So in what follows, we will assume that~$\delta$ is sufficiently small compared
to~$c-1$ for our estimates to hold.)
To show this, we will modify our argument from Section~\ref{sec:lowerdense2}.
Instead of partitioning the vertex set of~$G_{n,c/n}$ into two sets~$V'$ and~$V''$
and considering a long path in each of them, this time we will work with a single long path in
two stages. We will obtain this path in the first stage of a two-round exposure of~$G_{n,c/n}$.
More precisely, let $c_1:=(c+1)/2$. Thus $1<c_1<c$. Put $c_2:=(c-c_1)/(1-c_1/n)$.
So $c/n=c_1/n +c_2/n- c_1c_2/n^2$.
In our two-round exposure of~$G_{n,c/n}$ we first generate a random graph~$G_{n,c_1/n}$
and then independently~$G_{n,c_2/n}$. We then take their union, ignoring multiple
edges again. The choice of~$c_1$ ensures that there exists an $\alpha = \alpha (c)>0$
such that a.a.s.~$G_{n,c_1/n}$ contains a path~$P$ on~$2\alpha n$ vertices
(see~\cite{AKS}). We split~$P$ into two 
consecutive disjoint paths~$P'$ and~$P''$ each having~$\alpha n$ vertices.
We split~$P'$ further into $k':=2k=2\delta \sqrt{n}$ consecutive disjoint subpaths
$B_1,\ldots, B_{k'}$ each containing $t:=\alpha \sqrt{n}/(2 \delta)$ vertices.
The~$B_i$'s will be our candidate branch sets.

We will now use the greedy algorithm
described in Section~\ref{sec:greedy1} to join~$k$ of these candidate branch sets
into a $K_k$ minor in the second round of our two-round exposure of~$G_{n,c/n}$.
Thus similarly as before, let~$\U_0$ denote the set of all those pairs~$B_i,B_j$ of
candidate branch sets which are not joined by an edge in~$G_{n,c_2/n}$.
Let $U_0:=|\U_0|$ and let~$G'_0$ denote the auxiliary random graph whose vertices
are the~$B_i$ and whose edges correspond to the pairs in~$\U_0$.
So the edges of~$G'_0$ occur independently with probability $q:=(1-c_2/n)^{t^2}$.
Note that if~$\delta$ is sufficiently small compared to~$c-1$ then
$$ q=(1-c_2/n)^{\left(\alpha \sqrt{n}/(2\delta)\right)^2}\le
\exp(-c_2 \alpha^2/(4\delta^2))\le \delta^4
$$
and so
\begin{equation*}
\E (U_0)={{k'}^2 \choose 2}q \le 2\delta^2nq\le \delta^5n/2.
\end{equation*}
Again, we may assume that $\E(U_0)\ge n^{2/5}$. (In fact, an easy calculation shows
that~$\E(U_0)$ is linear in~$n$.) Thus
\begin{equation}\label{USize} 
\E(U_0)/2\le U_0 \le 2\E(U_0)\le \delta^5n
\end{equation}
by Chebyshev's inequality.
The path~$P''$ will play the role of the long path inside~$G_{n,p}[V'']$
considered in Section~\ref{sec:est}. We define~$i^*$ and $U_1,\dots,U_{i^*}$ as there and
divide~$P''$ into disjoint consecutive subpaths~$Q_1,\dots,Q_{i^*}$ where
\begin{equation}\label{eq:Qisparse}
|Q_i|=\frac{\delta^3 \alpha n}{2^{i}} \ \ \ \ \ \text{and so}\ \ \ \ \
|Q_i|=\frac{|Q_1|}{2^{i-1}}.
\end{equation}
We then split each~$Q_i$ into a set~$\mathcal{P}_i$ of smaller
paths, each containing~$\ell_i$ vertices, where~$\ell_i=|Q_i|/U_{i-1}$ as before.
So by~(\ref{USize}) we have that
\begin{equation}\label{eq:l1sparse}
\ell_1=\frac{|Q_1|}{U_0}\ge \frac{\alpha}{2\delta^2} \ \ \ \ \ \text{and}\ \ \ \ \
\ell_i=4^{i-1}\ell_1=\frac{4^{i-1}|Q_1|}{U_0}.
\end{equation}
Similarly as in Section~\ref{sec:greedy1}, our aim is to apply the greedy
algorithm for each~$i=1,\dots,i^*$ in turn to join all but precisely~$U_i$
of the pairs in~$\U_0$ that are still unjoined after the $(i-1)$th stage
of the algorithm, using one path in~$\mathcal{P}_i$ for each pair.
So suppose that $1\le i< i^*$ and after the $(i-1)$th stage of our algorithm we
have a set $\U_{i-1}\subseteq \U_0$ of unjoined pairs such that $|\U_{i-1}|=U_{i-1}$
and suppose that we now wish to run the $i$th stage of our algorithm. As before,
let~$G'_{i-1}$ denote the spanning subgraph of~$G'_0$ whose edge set corresponds
to~$\U_{i-1}$. Again, in the~$i$th stage of the algorithm we will first discard all those
pairs in~$\U_{i-1}$ which contain a candidate branch set lying in too many
pairs from~$\U_{i-1}$. This time, for all $i=1,\dots,i^*$ we put
\begin{equation}\label{eq:Deltaisparse}
\Delta_{i-1}= {8\ex(U_0) \over \delta^{1/2} 4^{i-1}k'}
\end{equation}
and let~$\B_{i-1}$ be the set of all those candidate branch sets whose degree
in~$G'_{i-1}$ is greater than~$\Delta_{i-1}$. We now run the algorithm described
in Section~\ref{sec:greedy1} and have to show that
$\prob(\geq U_{i-1}/8 \text{ paths in } \cP_i \text{ failed})$ is small. 

So as before, let~$\hat{p}$ denote the probability that a particular 
vertex is joined to a particular candidate branch set in~$G_{n,c_2/n}$.
Thus $\hat{p} = 1-(1-c_2/n)^t$. Since $1-c_2t/n\le (1-c_2/n)^t\le 1-c_2t/(2n)$
by~(\ref{eq:estp}) we have that
\begin{equation}\label{eq:hatpsparse}
\frac{c_2\alpha}{4\delta\sqrt{n}}=\frac{c_2 t}{2n}\le \hat{p}\le \frac{c_2 t}{n}
=\frac{c_2\alpha}{2\delta\sqrt{n}}.
\end{equation}
Then as before, the probability~$\hat{p}_i$ that a given path from~$\cP_i$ is joined to 
a particular candidate branch set satisfies~(\ref{eq:upperphat1})--(\ref{eq:lowerphat2}).
Thus
\begin{eqnarray*}
\hat{p}_i \Delta_{i-1} & \stackrel{(\ref{eq:upperphat1}), (\ref{eq:Deltaisparse})}{\le} &
\hat{p}\ell_i \frac{8\E (U_0)}{\delta^{1/2} 4^{i-1}k'}
\stackrel{(\ref{eq:hatpsparse}), (\ref{eq:l1sparse})}{\le}
\frac{c_2\alpha}{2\delta\sqrt{n}}\cdot  \frac{|Q_1|}{U_{0}}\cdot
\frac{4\E (U_0)}{\delta^{3/2}\sqrt{n}}
\stackrel{(\ref{eq:Qisparse}),(\ref{USize})}{\le}
\frac{2c_2\alpha^2 \delta^3 n}{\delta^{5/2}n}\le \delta^{1/3}
\end{eqnarray*}
if $\delta$ is sufficiently small compared to~$c-1$.
This gives an analogue to~(\ref{eq:hatpiCr}) and thus as before we can show
that~$(\ref{eq:failedP12})$ holds. Again, if $\hat{p}\ell_i \geq 1$, this
implies that $\prob(\geq U_{i-1}/8\text{ paths in }\cP_i \text{ failed}) \leq e^{-U_{i-1}}$.
So suppose that $\hat{p} \ell_i<1$ and thus $\hat{p}_i\ge \ell_i\hat{p}/2$
by~(\ref{eq:lowerphat1}). Hence
\begin{eqnarray*}
\hat{p}_i^2U_{i-1} \ge \left(\frac{\ell_i\hat{p}}{2}\right)^2 U_{i-1}
\stackrel{(\ref{eq:l1sparse})}{=}\frac{4^{i-1} \ell_1|Q_i|}{4}\hat{p}^2
\stackrel{(\ref{eq:hatpsparse})}{\ge}
\frac{4^{i-1} \ell_1 |Q_i|}{4 }\cdot \frac{c_2^2\alpha^2}{16\delta^2 n}
\stackrel{(\ref{eq:Qisparse})}{\ge} \frac{c_2^2\delta \alpha^3\ell_1}{2^7}
\stackrel{(\ref{eq:l1sparse})}{\ge} \frac{c_2^2 \alpha^4}{2^{8}\delta}\ge 2^{8}
\end{eqnarray*}
if $\delta$ is sufficiently small compared to~$c-1$.
So~(\ref{eq:hatpiUi}) still holds and as before this shows that
$$\prob(\geq U_{i-1}/8\text{ paths in }\cP_i \text{ failed}) \leq e^{-U_{i-1}}
=o(1/\ln n)$$
holds in this case too. As before this implies that a.a.s. our algorithm
produces sets~$\U_i$ with $|\U_i|=U_i$ for all $i=1,\dots,i^*$.
Similarly as in~(\ref{eq:sumBi}) one can show that $\sum_{j\ge 0} |\B_j|\le \delta^{1/2}k'$.
Thus by deleting all the candidate branch sets in $\bigcup_{j\ge 0} \B_j$
as well as one candidate branch set from every pair in~$\U_{i^*}$
we obtain a complete graph of order~$k'-U_{i^*}- |\bigcup_{j\ge 0} \B_j|
\ge k$ as a minor. This completes the proof of the lower bound in Theorem~\ref{thm2}.

\medskip

{\footnotesize \obeylines \parindent=0pt

Nikolaos Fountoulakis, Daniela K\"{u}hn \& Deryk Osthus 
School of Mathematics
University of Birmingham
Edgbaston
Birmingham
B15 2TT
UK
}

{\footnotesize \parindent=0pt

\it{E-mail addresses}:
\tt{\{nikolaos,kuehn,osthus\}@maths.bham.ac.uk}}

\end{document}